\documentclass[preprint, 3p, 11pt]{elsarticle}

\usepackage{amssymb,amsmath}
\usepackage{graphicx,color,epsfig}

\newenvironment{prf}[1][Proof]{\noindent\textbf{#1.} }{\hfill$\square$}

\newcommand{\R}{\mathbb{R}}

\newcommand{\Rn}{\mathbb{R}^{n}}
\newcommand{\Rm}{\mathbb{R}^{m}}
\newcommand{\Rd}{\mathbb{R}^{d}}

\newcommand{\Rnd}{\mathbb{R}^{n\times d}}

\newcommand{\PP}{\mathbb{P}}

\newcommand{\E}{\mathbb{E}}

\newcommand{\F}{\mathcal{F}}

\newcommand{\U}{\mathcal{U}}

\newcommand{\tr}{\mbox{tr}}

\renewcommand{\le}{\leqslant}
\renewcommand{\ge}{\geqslant}
\renewcommand{\leq}{\leqslant}
\renewcommand{\geq}{\geqslant}

\newtheorem{thm}{Theorem}[section]
\newtheorem{prop}[thm]{Proposition}

\newtheorem{defn}[thm]{Definition}
\newtheorem{cor}[thm]{Corollary}

\makeatletter
\def\ps@pprintTitle{%
  \let\@oddhead\@empty
  \let\@evenhead\@empty
  \def\@oddfoot{\reset@font\hfil\thepage\hfil}
  \let\@evenfoot\@oddfoot
}
\makeatother

\begin{document}

\begin{frontmatter}
\title{On the Stability of Receding Horizon Control for\\ Continuous-Time Stochastic Systems}

\author[les]{Fajin Wei}
\ead{F.Wei@alumni.lboro.ac.uk}
\author[les]{Andrea Lecchini-Visintini\corref{cor}}
\ead{alv1@leicester.ac.uk}
\address[les]{Department of Engineering, University of Leicester, Leicester, LE1 7RH, UK.}
\cortext[cor]{Corresponding author}

\begin{abstract}
We study the stability of receding horizon control for continuous-time non-linear stochastic differential equations.
We illustrate the results with a simulation example in which
we employ receding horizon control to design an investment strategy to repay a debt.
\end{abstract}

\begin{keyword}
Receding horizon control, Stochastic differential equations, Stochastic optimal control, Hamilton-Jacobi-Bellman equations,
Lyapunov functions, It\^{o}'s formula, optimal investment.
\end{keyword}

\end{frontmatter}

\section{Introduction}
\label{intro}
In Receding Horizon Control (RHC), the control action, at each time $t$ in $[0,\infty)$,
is derived from the solution of an optimal control
problem defined over a finite future horizon $[t,\, t+T]$.
The RHC strategy establishes a feedback law which, under
certain conditions, can ensure asymptotic stability of the controlled system.
This control strategy has been successfully developed over
the last twenty years for systems described by deterministic equations.
In this context RHC is also well known as Model Predictive Control (MPC)
and has proven to be very successful in dealing with non-linear and constrained
systems, see e.g.~\cite{Mayne-et-al-00,Maciejowski-02,Magni-Scattolini-04}.
The extension of RHC from deterministic to stochastic systems is the objective of current research.
RHC schemes for the control of discrete-time
stochastic systems have been proposed recently in
\cite{Primbs-et-Sung-09,Cannon-et-al-09,Chatterjee-et-al-11,Hokayem-et-al-12}.

In this note, we discuss RHC for systems described by  continuous-time non-linear
stochastic differential equations (SDEs). To the extent of our knowledge,
the RHC strategy has not yet been considered in this context.
In order to study the stability of RHC for continuous-time SDEs,
we formulate conditions under which the value function of the associated
finite-time optimal control problem can be used as Lyapunov function
for the RHC scheme.
This is a well established approach for studying the stability
of RHC schemes, which here is extended using Lyapunov criteria for
stochastic dynamical systems \cite{Kushner-65}.
We illustrate this contribution with a simple example of
an optimal investment problem.
Optimal investments problems are well suited to  be tackled by stochastic
control methods, see e.g.~\cite{Primbs-09,Pola-et-Pola-09}.
In our example, we design an investment strategy to repay a debt.
Having negative wealth due to an initial debt, the investor has the option
to increase his/her current debt in order to buy a risky asset.
The asymptotic stability of the adopted RHC scheme
guarantees that the wealth of the investor tends to zero,
so that the initial debt is eventually repaid.


\section{Problem statement}

Let $(\Omega,\F,\PP)$ be a complete probability space equipped with the natural filtration $(\F_{t})_{t\ge0}$
generated by a standard Wiener process $W:[0,\infty)\times\Omega\to\Rd$ on it.
We consider a controlled time-homogeneous SDE for a process $X:[0,\infty)\times\Omega\to\Rn$,
\begin{eqnarray}\label{eq:sde1}
dX^{0,x_{0},u}_{t}&=&b(X^{0,x_{0},u}_{t},u_{t})dt+\sigma(X^{0,x_{0},u}_{t},u_{t})dW_{t},\\
X^{0,x_{0},u}_{0}&=&x_{0},\nonumber
\end{eqnarray}
where $x_{0}\in\Rn$; $b:\Rn\times\Rm\to\Rn$ and $\sigma:\Rn\times\Rm\to\Rnd$
are measurable functions and satisfy
$$
|b(x,u)|+|\sigma(x,u)|\le C(1+|x|),\forall(x,u)\in\Rn\times U,\mbox{(linear growth)},
$$
and
$$
|b(x,u)-b(y,u)|+|\sigma(x,u)-\sigma(y,u)|\le C|x-y|,\forall(x,y,u)\in\Rn\times\Rn\times U,\mbox{(Lipschitz)},
$$
for some constant $C$; and $u_{(\cdot)}$ is an admissible control process
$$u_{(\cdot)}\in\U:=\left\{u:[0,\infty)\times\Omega\to U:\mbox{ progressively measurable and }\E\int_{0}^{\infty}|u_{t}(\omega)|^{2}dt<\infty\right\},$$
with the set $U\subset\Rm$ compact.
Here the superscripts of $X^{0,x_{0},u}$ mean that the initial value of the process at time $0$ is $x_{0}$
and the involved control process is $u_{(\cdot)}$.
In this paper we are concerned with the conditions under which there exists a control process that drives the stochastic system $X$
to the origin $0\in\Rn$ and guarantees asymptotic stability of the controlled process.
Here, the following definition of stability is adopted \cite{Kushner-65}:
\begin{defn}\label{def:stab}\upshape
Given a stochastic continuous-time process $X: \R_{+}\times\Omega\to\Rn$, where $\R_{+}:=[0,\infty)$, with $X_{0}=x_{0}\in\Rn$
\begin{itemize}
\item[(S1)] The origin is stable almost surely
if and only if, for any $\rho>0$, $\epsilon>0$, there is a $\delta(\rho,\epsilon)>0$ such that,
if $|x_{0}|\le\delta(\rho,\epsilon)$,
    $$
    \PP\left[\sup_{t\in\R_{+}}|X_{t}|\ge\epsilon\right]\le\rho.
    $$
\item[(S1')]An equivalent definition to (S1) is:
Let $h(\cdot):\R_{+}\to\R_{+}$ be a scalar-valued, nondecreasing,
and continuous function of $|x|$. Let $h(0)=0$, $h(r)>0$ for $r\ne0$.
Then the origin is stable almost surely if and only if, for any
$\rho>0$, $\lambda>0$, there is a $\delta(\rho,\lambda)>0$ such that,
for $|x_{0}|\le\delta(\rho,\lambda)$,
    $$
    \PP\left[\sup_{t\in\R_{+}}h(|X_{t}|)\ge\lambda\right]\le\rho.
    $$
\item[(S2)] The origin is asymptotically stable almost surely if and only if
it is stable a.s., and $X_{t}\to0$ a.s. for all $x_{0}$ in some neighborhood $R$ of the origin.
If $R=\Rn$ then we add `in the large'.
\end{itemize}
\end{defn}


\section{Main results}\label{sec:main}

Let $T>0$. As a preliminary step we consider the SDE
for $X:[t,T]\times\Omega\to\Rn$ starting from the point
$x\in\Rn$ at the time $t\in[0,T]$
\begin{eqnarray}\label{eq:sde2}
dX^{t,x,u}_{s}&=&b(X^{t,x,u}_{s},u_{s})ds+\sigma(X^{t,x,u}_{s},u_{s})dW_{s},\\
X^{t,x,u}_{t}&=&x.\nonumber
\end{eqnarray}
Let $f:\Rn\times\Rm\to\R_{+}$ and $g:\Rn\to\R_{+}$ be measurable nonnegative functions.
Now we consider the problem of minimizing the following cost functional, $\forall(t,x)\in[0,T]\times\Rn$,
\begin{equation}\label{eq:cost}
J[t,x;T;u_{(\cdot)}]:=\E\left[\int_{t}^{T}f(X^{t,x,u}_{s},u_{s})ds+g(X^{t,x,u}_{T})\right]
\end{equation}
over the set $\U$ of admissible control processes.
We define the value function as
\begin{equation}\label{eq:valuefunc}
v(t,x;T):=\inf_{u_{(\cdot)}\in\U}J[t,x;T;u_{\cdot}]=\inf_{u_{(\cdot)}\in\U}\E\left[\int_{t}^{T}f(X^{t,x,u}_{s},u_{s})ds+g(X^{t,x,u}_{T})\right].
\end{equation}
and denote $u^*_s(t,x;T)$, $t\leq s\leq T$, the optimal control process if it exists.
In particular, when $t=0$ we denote $V(x;T):=v(0,x;T)$.

Standard stochastic optimal control theories (see, for instance, \cite{Fleming-et-Rishel-75,Fleming-et-Soner-06})
about the controlled SDE \eqref{eq:sde1} tell us that
the Hamilton-Jacobi-Bellman (HJB) equation for the value function
$v(\cdot,\cdot;T)$ is, $\forall(t,x)\in[0,T]\times\Rn$,
\begin{eqnarray}\label{eq:hjb1}
-\partial_{t}v(t,x;T)&=&\inf_{u\in U}\left[\frac{1}{2}\tr[\sigma\sigma^{*}(x,u)D^{2}v(t,x;T)]+\langle b(x,u),Dv(t,x;T)\rangle+f(x,u)\right],\\
v(T,x;T)&=&g(x).\nonumber
\end{eqnarray}
Hereafter we use the notations
$$\partial_{t}v:=\frac{\partial v}{\partial t},\
Dv=\left(\begin{array}{c}
\frac{\partial v}{\partial x_{1}}\\\vdots\\\frac{\partial v}{\partial x_{n}}
\end{array}\right),
\mbox{ and }
D^{2}v=
\left(\begin{array}{cccc}
\frac{\partial^{2}v}{\partial x_{1}^{2}}&\frac{\partial^{2}v}{\partial x_{1}\partial x_{2}}
&\cdots&
\frac{\partial^{2}v}{\partial x_{1}\partial x_{n}}\\
\vdots&\vdots&\cdots&\vdots\\
\frac{\partial^{2}v}{\partial x_{n}\partial x_{1}}&\frac{\partial^{2}v}{\partial x_{n}\partial x_{2}}
&\cdots&
\frac{\partial^{2}v}{\partial x_{n}^{2}}
\end{array}\right).
$$
Suppose this HJB equation has a solution and that the infimum in the equation is attained by $\tilde{u}(t,x;T)$
for every $(t,x)\in[0,T]\times\Rn$, i.e.,
$$
-\partial_{t}v(t,x;T)=\frac{1}{2}\tr[\sigma\sigma^{*}(x,\tilde{u}(t,x;T))D^{2}v(t,x;T)]+\langle b(x,\tilde{u}(t,x;T)),Dv(t,x;T)\rangle+f(x,\tilde{u}(t,x;T)),
$$
then we construct the optimal control process for the SDE \eqref{eq:sde2} as
\begin{equation}\label{eq:ustar}
u^{*}_{s}(t,x;T):=\tilde{u}(s,X^{t,x,\tilde{u}}_{s};T),\ t\le s\le T.
\end{equation}
As a consequence, the value function turns out to be
$$
v(t,x;T)=\E\left[\int_{t}^{T}f(X^{t,x,u^{*}}_{s},u^{*}_{s})ds+g(X^{t,x,u^{*}}_{T})\right]=\E\left[\int_{t}^{T}f(X^{t,x,\tilde{u}}_{s},\tilde{u}(s,X^{t,x,\tilde{u}}_{s};T))ds+g(X^{t,x,\tilde{u}}_{T})\right].
$$
Now, for all the states $x\in\Rn$ all the time, we apply the specifically designed feedback law
\begin{equation}\label{eq:RHC}
u^{c}(x;T):=\tilde{u}(0,x;T)
\end{equation}
to the stochastic system \eqref{eq:sde1}.
In other words, for the state $X^{0,x_{0},u^{c}}_{t}$, at any time $t\ge0$,
we apply only the initial optimal control
$$
u^{c}(X^{0,x_{0},u^{c}}_{t};T)=\tilde{u}(0,X^{0,x_{0},u^{c}}_{t};T)=u^{*}_{0}(0,X^{0,x_{0},u^{c}}_{t};T)
$$
to the system.
In particular, for $t=0$,
$u^c(X^{0,x_{0},u^{c}}_{0};T)=u^{c}(x_{0};T)=\tilde{u}(0,x_{0};T)$.
We call $u^{c}(\cdot;T)$ the
continuous receding horizon control
process with the receding horizon $T$ \mbox{for the controlled SDE \eqref{eq:sde1}.}

When $t=0$, using the HJB equation and \eqref{eq:RHC}, we obtain
$$
-\partial_{t}v(t,x;T)|_{t=0}-f(x,u^{c}(x;T))=\frac{1}{2}\tr[\sigma\sigma^{*}(x,u^{c}(x;T))D^{2}V(x;T)]+\langle b(x,u^{c}(x;T)),DV(x;T)\rangle.
$$
Let us denote
\begin{equation}\label{eq:phi}
\phi(x;T):=\partial_{t}v(t,x;T)|_{t=0}+f(x,u^{c}(x;T))
\end{equation}
for all $x\in\Rn$.
Let us assume that
\begin{equation}\label{A1}
V\in C^{2,1}(\Rn\times\R_{+},\R),\tag{A1}
\end{equation}
i.e.~the set of functions $\psi:\Rn\times\R_{+}\to\R,\ (x,t)\mapsto\psi(x,t)$
that are twice continuously differentiable in $x$ and once continuously differentiable in $t$; and that
\begin{equation}\label{A2}
\phi(x;T)>0,\ \forall x\in\Rn\setminus\{0\},\mbox{ and }\phi(0;T)=0.\tag{A2}
\end{equation}
In particular, note that $\phi(0;T)=0$ can be achieved if we suppose that
\begin{equation}\label{A2.1}
b(0,0)=0, \mbox{ and }  \sigma(0,0)=0;\tag{A2.1}
\end{equation}
and that
\begin{equation}\label{A2.2}
f(0,0)=0,\,   f(x,u)>0,\, \forall x\in\Rn\setminus\{0\},\ \forall u\in U,
\,\mbox{ and }\,
g(0)=0,\, g(x)>0,\,\forall x\in\Rn\setminus\{0\}.\tag{A2.2}
\end{equation}
We now state one of our main results and then discuss Assumptions (\ref{A1}-\ref{A2})
in more detail.
\begin{prop}[Convergence]\label{th:converge}\upshape
Suppose $b(\cdot,\cdot)$ and $\sigma(\cdot,\cdot)$ are of linear growth and Lipschitz.
Suppose the HJB equation \eqref{eq:hjb1} has a unique classical solution.
Under the assumptions \eqref{A1}, and \eqref{A2},
we have that almost all the trajectories of the stochastic system \eqref{eq:sde1}
driven by the continuous receding horizon control $u^{c}(\cdot;T)$ defined in \eqref{eq:RHC}
converge to the origin.
\end{prop}

\begin{prf}
Temporarily we denote $X_{t}:=X^{0,x_{0},u^{c}}_{t}$ for simplicity.
By It\^{o}'s formula and the Hamilton-Jacobi-Bellman equation we have
\begin{eqnarray*}
dV(X_{t};T)&=&\langle DV(X_{t};T),dX_{t}\rangle+\frac{1}{2}\tr[\sigma\sigma^{*}(X_{t},u^{c}(X_{t};T))D^{2}V(X_{t};T)]dt\\
&=&\langle DV(X_{t};T),b(X_{t},u^{c}(X_{t};T))\rangle dt+\langle DV(X_{t};T),\sigma(X_{t},u^{c}(X_{t};T))dW_{t}\rangle\\
&& +\,\,\frac{1}{2}\tr[\sigma\sigma^{*}(X_{t},u^{c}(X_{t};T))D^{2}V(X_{t};T)]dt\\
&=&-\phi(X_{t};T)dt+\langle DV(X_{t};T),\sigma(X_{t},u^{c}(X_{t};T))dW_{t}\rangle.
\end{eqnarray*}
Therefore by the assumption \eqref{A2} we get that, for $0\le s\le t<\infty$,
\begin{equation}\label{V-supermart}
\E[V(X_{t};T)|\F_{s}]-V(X_{s};T)=-\E\left[\int_{s}^{t}\phi(X_{\xi};T)d\xi\Big|\F_{s}\right]\le0.
\end{equation}
In particular when $s=0$, $\E \left[V(X_{t};T)\right]\le V(X_{0};T)=V(x_{0};T)<\infty$.
Hence, $\{V(X_{t};T),\F_{t}\}_{0\le t<\infty}$ is a nonnegative supermartingale.
Let us recall the supermartingale convergence theorem, see, for instance, \cite[pag.~18]{Karatzas-et-Shreve-06}:
Suppose $\{M_{t},\F_{t}\}_{0\le t<\infty}$ is a right-continuous, nonnegative supermartingale.
Then $M_{\infty}(\omega):=\lim_{t\to\infty}M_{t}(\omega)$ exists for $\PP$-a.e. $\omega\in\Omega$,
and $\{M_{t},\F_{t}\}_{0\le t\le\infty}$ is a supermartingale.
Hence, there exists a random variable $V_{\infty}$, integrable, such that
$V(X_{t};T)\to V_{\infty},\mbox{ a.s.}$.
Now, equation \eqref{V-supermart} implies particularly that
\begin{equation}
\E \left[V_{\infty}\right]-V(x_{0};T)=-\int_{0}^{\infty}\E\left[\phi(X_{t};T)\right]dt,
\end{equation}
which in turn gives
$\lim_{t\to\infty}\E\left[\phi(X_{t};T)\right]=0$.
Therefore we obtain that
$$
\lim_{t\to\infty}\phi(X_{t};T)=0,\mbox{ a.s.,}
$$
because $\phi\ge0$. Denote the set $\Phi:=\{x\in\Rn:\phi(x;T)=0\}$ and from the assumption \eqref{A2}
we know that $\Phi=\{0\}$ and it is closed. Thus by the continuity of $\phi$ we learn that almost all
the trajectories of the process $X_{t}$ converge to the origin in $\Rn$.\smallskip
\end{prf}

In \eqref{A1}, continuous twice differentiability of $V$ is required for
the applicability of It\^{o}'s formula.
In order to illustrate Assumption (\ref{A2}) note that, since \eqref{eq:sde2}
is time-homogeneous, we have
$$
\partial_{t}v(t,x;T)|_{t=0} =  - \partial_{H}v(0,x;H)|_{H=T} =  - \partial_{H} V(x;H)|_{H=T}.
$$
Hence, $\phi(x;T)$ can be equivalently expressed as:
\begin{equation}
\phi(x;T) = - \partial_{H} V(x;H)|_{H=T} + f(x,u^{c}(x;T)).
\end{equation}
Note that, here, $\partial_{H} V(x;H)|_{H=T}$ is the rate at which the optimal cost
increases with increments in the horizon $T$.
Note that $\phi(x;T)>0$ is implied by
\begin{equation}\label{eq:neg}
\partial_{H} V(x;H)|_{H=T} \leq 0\, .
\end{equation}
Here, we recover the condition of monotonic decrease of
the value function with the length of the horizon, which is
a well-studied condition for the stability of receding horizon
control schemes in a deterministic setting, see e.g.~\cite{Magni-Scattolini-04,Chen-Shaw-82,Bertsekas-05}.
Finally,  note that Assumption (\ref{A2.1}) requires that the
state, at the origin, is not affected by the noise.
This is an unavoidable  assumption  for obtaining
asymptotic stability in the sense of Definition \ref{def:stab}.
If this assumption is not met, then one has to resort to other notions
of stability such as, for example,\mbox{ mean square boundeness,
see e.g.~\cite{Chatterjee-et-al-11,Hokayem-et-al-12,Chatterjee-et-al-12}.}

Under additional assumptions, asymptotic stability according to Definition \ref{def:stab} can be obtained.
Let us assume in addition that for any $\epsilon>0$, there exists a $\delta(\epsilon)>0$ such that
\begin{equation}\label{A3}
V(x;T)<\epsilon \quad\mbox{if } |x|\le\delta(\epsilon); \tag{A3}
\end{equation}
and that there exists
a continuous, nondecreasing function $h:\R_{+}\to\R_{+}$, satisfying
$$
h(0)=0,\ h(r)>0, \forall r\ne0,
$$
such that
\begin{equation}\label{A4}
V(x;T)\ge h(|x|),\ \forall x\in\Rn.\tag{A4}
\end{equation}
\begin{cor}[Asymptotic stability]\upshape\label{th:corr}
Suppose $b(\cdot,\cdot)$ and $\sigma(\cdot,\cdot)$ are of linear growth and Lipschitz.
Under the assumptions \eqref{A1}, \eqref{A2}, \eqref{A3}, and \eqref{A4},
we have that the origin in $\Rn$ is asymptotically stable.
\end{cor}
\begin{prf}
Let us recall the supermartingale inequality, see, \cite[pag.~13]{Karatzas-et-Shreve-06}:
Suppose $\{M_{t},\F_{t}\}_{0\le t<\infty}$ is a supermartingale whose every path is right-continuous.
Let $\lambda>0$ and $[a,b]$ be a subinterval of $[0,\infty)$. Then, we have
$$
\PP\left[\sup_{a\le t\le b}M_{t}\ge\lambda\right]\le\frac{\E[M^{+}_{a}]}{\lambda}\, .
$$
Hence, for any $\lambda>0$, we have
$$
\PP\left[\sup_{0\le t<\infty}V(X_{t};T)\ge\lambda\right]\le\frac{\E[V(x_{0};T)]}{\lambda}=\frac{V(x_{0};T)}{\lambda}\, .
$$
Using assumption \eqref{A3}, we obtain that for any $\lambda>0$, $\rho>0$, there
exists $\delta(\rho,\lambda)>0$ such that
$$
\PP\left[\sup_{0\le t<\infty}V(X_{t};T)\ge\lambda\right]\le\rho \quad\mbox{if } |x_{0}|<\delta(\rho,\lambda).
$$
Then, by assumption \eqref{A4}, we immediately obtain
$$\PP\left[\sup_{0\le t<\infty} h(X_{t})\ge\lambda\right]\le\PP\left[\sup_{0\le t<\infty}V(X_{t};T)\ge\lambda\right]\le\rho,$$
which entails the asymptotic stability of the origin in $\Rn$.
\end{prf}

In the following section, we present a simple illustrative example in which
Assumptions  (\ref{A1}-\ref{A4})  can be verified explicitly.


\section{Example: repayment of a debt}\label{sec:example}

In this example we consider a variant of the
Merton's portfolio problem, see e.g.~\cite{Fleming-et-Rishel-75,Fleming-et-Soner-06}.
Suppose in a complete financial market there are only two assets, one asset being risk free such as, for instance,
a bank deposit or a bond, and the other one being a risky asset such as, for instance, a stock.
The assets obey the following price process, $t\ge0$,
\begin{eqnarray}\label{eq:exprices}
dP^{1}_{t}&=&rP^{1}_{t}dt,\\
dP^{2}_{t}&=&bP^{2}_{t}dt+\sigma P^{2}_{t}dW_{t},\nonumber
\end{eqnarray}
where $b>r>0$, and $\sigma\ne0$ are given constants.  Here $r$ is the risk free interest rate for the bank deposit,
$b$ is the drift, or average, rate of the stock's return, and $\sigma$ is  the volatility of the stock's return.
Suppose an agent's total wealth at time $t$ is $X_{t}$, comprising the risk
free part $\Pi^{1}_{t}P^{1}_{t}$ and the risky part $\Pi^{2}_{t}P^{2}_{t}$,
where $\Pi^{1}_{t}$ and $\Pi^{2}_{t}$ are the quantity of the assets, respectively.
The control variable $u_t$ is the portion of the total wealth $X_{t}$ invested by the agent on
the risky asset at time $t$. In a continuous-time setting, it is assumed that the allocation
of wealth takes place instantaneously. Hence, for $t\ge0$, we have:
\begin{eqnarray}\label{eq:exassets}
X_{t}&=&\Pi^{1}_{t}P^{1}_{t}+\Pi^{2}_{t}P^{2}_{t},\nonumber \\
\Pi^{1}_{t}P^{1}_{t}&=&(1-u_{t})X_{t},\\
\Pi^{2}_{t}P^{2}_{t}&=& u_{t}X_{t}.\nonumber
\end{eqnarray}
We assume that the investment obeys the self financing condition.
That is, starting from an initial wealth $x_{0}$, the agent can only sell or buy these
two assets but is not allowed to borrow money from outside or consume his or her wealth.
Written down as a differential equation, this is $P^{1}_{t}d\Pi^{1}_{t}+P^{2}_{t}d\Pi^{2}_{t}=0$.
Using \eqref{eq:exprices} and \eqref{eq:exassets} we obtain
\begin{eqnarray*}
dX_{t}&=&\Pi^{1}_{t}dP^{1}_{t}+\Pi^{2}_{t}dP^{2}_{t}+P^{1}_{t}d\Pi^{1}_{t}+P^{2}_{t}d\Pi^{2}_{t}\\
&=&\Pi^{1}_{t}dP^{1}_{t}+\Pi^{2}_{t}dP^{2}_{t}\\
&=&[r+(b-r)u_{t}]X_{t}dt+\sigma u_{t}X_{t}dW_{t}.
\end{eqnarray*}
Therefore, the wealth satisfies the stochastic differential equation
\begin{eqnarray}\label{eq:exwealth}
dX^{0,x_{0},u}_{t}&=&[r+(b-r)u_{t}]X^{0,x_{0},u}_{t}dt+\sigma u_{t}X^{0,x_{0},u}_{t}dW_{t},\\
X^{0,x_{0},u}_{0}&=&x_{0},\nonumber
\end{eqnarray}
which is called the wealth process.

In this example, we assume that the agent's initial wealth is negative and his or
her aim is to repay any debt eventually.
Hence, we have $x_{0}<0$ and investigate the asymptotic stability of the origin $0\in\R$.
In this case, the wealth process will always be negative until the time it reaches zero.
However, note that $X_{t}<0$ does not necessarily mean that the
agent has only debt without any money to invest or stock to sell.
If negative values of $u_{t}$ in \eqref{eq:exassets} are allowed then
the investor is able to increase his or her debt in order to buy stocks.
Let us explain the financial meanings of all possible values
of $u_{t}$ when $X_{t}<0$:
\begin{itemize}
\item If $u_t=0$, then  $\Pi^{1}_{t}P^{1}_{t}= X_{t}$ and $\Pi^{2}_{t}P^{2}_{t}= 0$.
In this case, the investor just owns a debt with the bank equal to his or her negative wealth.
\item If $u_t < 0$, then $\Pi^{1}_{t}P^{1}_{t} < X_{t}$ and $\Pi^{2}_{t}P^{2}_{t} > 0$.
In this case, the investor is borrowing additional funds from the bank and is using
it to buy stocks.
The investor owns a debt with the bank equal to $(1-u_t)X_t$ and
a positive quantity of stocks whose value is $u_t X_{t}$.
\item If $u_t > 0$, then $\Pi^{1}_{t}P^{1}_{t} > X_{t}$ and $\Pi^{2}_{t}P^{2}_{t} < 0$.
In this case, the investor is borrowing stocks (short selling) and is using this additional
wealth to reduce his or her debt with the bank. However, in general, this is not desirable
because a debt in stocks is more risky than a debt with the bank.
\end{itemize}
Note that if $\Pi^{1}_{t}P^{1}_{t} <0$ (i.e. when the investor has a debt with the bank)
then $r$ is the rate at which interests are payed to the  the bank when owing debt.
Here, in order to simplify the exposition of the problem, we assume that $r$ is the same
whether $\Pi^{1}_{t}P^{1}_{t} <0$ or $\Pi^{1}_{t}P^{1}_{t} >0$.
However, it will be shown that this issue is immaterial. In fact, the derived control process
will be constantly negative until the wealth reaches zero. In turn, this means that
$\Pi^{1}_{t}P^{1}_{t}<0$ and, therefore, $r$ will not change meaning throughout.
Finally, we  consider the constraints $u_{t}\in[c_{1},c_{2}]$ with $c_{1}<0$ and $c_{2}\geq 0$.
The constraint $c_{1}$ means that the investor is not allowed to increase the
debt with the bank by more than $(1+|c_{1}|)$ times his or her current (negative) wealth.
Similarly, $c_{2}$ is a constrain on borrowing stocks. If short selling is not allowed
then $c_{2}=0$.

\subsection{Running cost}\label{sec:exrun}

We consider  cost function \eqref{eq:cost} with
\begin{equation}\label{eq:excost}
f(x,u):=\begin{cases}
(-x)^{\beta},&x\le0,\\
0,&x>0,
\end{cases}\quad\mbox{ and }\quad
g(x):=0,\,\forall x,
\end{equation}
with $\beta>2$ being a given constant.
In this case, the value function $v(t,x;T)$
satisfies the HJB equation, $(t,x)\in[0,T]\times(-\infty,0]$,
\begin{eqnarray}\label{eq:ex1HJB}
-\partial_{t}v(t,x;T)&=&\inf_{u\in[c_{1},c_{2}]}\left[\frac{1}{2}(\sigma ux)^{2}D^{2}v(t,x;T)+(r+(b-r)u)xDv(t,x;T)+(-x)^{\beta}\right],\\
v(T,x;T)&=&0.\nonumber
\end{eqnarray}
If $D^{2}v(t,x;T)>0$ then a necessary condition for $\tilde{u}$ to be a minimiser is
$$
\sigma^{2}\tilde{u}x^{2}D^{2}v(t,x;T)+(b-r)xDv(t,x;T)=0,
$$
that is,
\begin{equation}\label{eq:ex1utilde}
\tilde{u}=-\frac{(b-r)Dv(t,x;T)}{\sigma^{2}xD^{2}v(t,x;T)}.
\end{equation}
In order to solve the HJB equation \eqref{eq:ex1HJB}, we try to find a value function in the form
\begin{equation}\label{eq:ex1w}
v(t,x;T)=\begin{cases}
(-x)^{\beta}w(t),&x\le0,\\
0,&x>0,
\end{cases}
\end{equation}
with $w$ to be determined.
For $x\le0$, by substituting \eqref{eq:ex1w} into \eqref{eq:ex1utilde},
we obtain
\begin{equation}\label{eq:utilde}
\tilde{u}=-\frac{(b-r)(-\beta)(-x)^{\beta-1}}{\sigma^{2}x\beta(\beta-1)(-x)^{\beta-2}}=-\frac{b-r}{(\beta-1)\sigma^{2}}.
\end{equation}
Note that the so-obtained $\tilde{u}$ has the same expression as
in the classical Merton's problem (although here we assumed $\beta>2$ instead of $\beta<1$),
see e.g.  \cite[pag.~160-161]{Fleming-et-Rishel-75}, or \cite[pag.~168-169]{Fleming-et-Soner-06}.
Using \eqref{eq:ex1w} and \eqref{eq:utilde},
the HJB equation becomes
\begin{equation*}
-(-x)^{\beta}w'(t)=\frac{\beta(b-r)^{2}}{2(\beta-1)\sigma^{2}}x^{2}(-x)^{\beta-2}w(t)-\left(\beta r-\frac{\beta(b-r)^{2}}{(\beta-1)\sigma^{2}}\right)(-x)^{\beta-1}xw(t)+(-x)^{\beta}
\end{equation*}
that is,
\begin{equation}\label{ex1:wdiff}
w'(t)=\left[\frac{\beta(b-r)^{2}}{2(\beta-1)\sigma^{2}}-\beta r\right]w(t)-1
\end{equation}
and we obtain that the solution is
\begin{equation}\label{eq:ex1wsol}
w(t)=\frac{1}{\eta}\left(1-e^{\eta(t-T)}\right),
\end{equation}
where
\begin{equation}\label{eq:ex1eta}
\eta:=\frac{\beta(b-r)^{2}}{2(\beta-1)\sigma^{2}}-\beta r\, .
\end{equation}
Hence, the value function \eqref{eq:ex1w} is actually given by
\begin{eqnarray}\label{eq:ex1value}
v(t,x;T)=\begin{cases}
(-x)^{\beta}\frac{1}{\eta}\left(1-e^{\eta(t-T)}\right),&x\le0,\\
0,&x>0.
\end{cases}
\end{eqnarray}
Here we assume $\eta\ne0$.
Note that whether $\eta<0$ or $\eta>0$ it always holds that $w(t)>0$.
Thus $D^{2}v=\beta(\beta-1)(-x)^{\beta-2}w(t)>0$ as we expected.
(However, in the next subsection, we will narrow the requirement to be $\eta>0$.)

Eventually, we obtain that the corresponding receding horizon
control process is
\begin{equation}\label{eq:ex1uc}
u^{c}(X_{t};T)=\tilde{u}(0,X_{t};T)=-\frac{b-r}{(\beta-1)\sigma^{2}}.
\end{equation}
The stochastic system \eqref{eq:exwealth} under the receding horizon control process
\eqref{eq:ex1uc} becomes
\begin{eqnarray}
dX^{0,x_{0},u}_{t}&=&\left[r-\frac{(b-r)^{2}}{(\beta-1)\sigma^{2}}\right]X^{0,x_{0},u}_{t}dt-\frac{b-r}{(\beta-1)\sigma}X^{0,x_{0},u}_{t}dW_{t},\\
X^{0,x_{0},u}_{0}&=&x_{0}.\nonumber
\end{eqnarray}
In this case the solution turns out to be a geometric Brownian motion that can be
written down explicitly (see e.g.\cite[pag.~349-50]{Karatzas-et-Shreve-06})
\begin{eqnarray}\label{eq:ex1sol}
X^{0,x_{0},u}_{t}&=&
x_{0}\exp\left\{\left[r-\frac{(2\beta-1)(b-r)^{2}}{2(\beta-1)^{2}\sigma^{2}}\right]t-\frac{b-r}{(\beta-1)\sigma}W_{t}\right\}.
\end{eqnarray}
The asymptotic properties near the origin can be seen directly from the explicit solution \eqref{eq:ex1sol}.
However, in the following subsection we will use Proposition \ref{th:converge}  and Corollary \ref{th:corr}
to asses the stability of the system.
Note that the derived control policy \eqref{eq:ex1uc} is in fact a constant process with negative value.
This implies that the agent is always advised to borrow money from the bank and buy stocks with it.
He or she will repay the debt in the end when his or her total wealth reaches zero by making profit
from investment in stocks.

\subsection{Verification of assumptions}

Here we verify when Assumptions (\ref{A1}-\ref{A4}) are met by the receding horizon control process  \eqref{eq:ex1uc}.
Note that, for $v(0,x;T)$ given by $\eqref{eq:ex1value}$ we have $V(x;T)=v(0,x;T)<\infty$ for all $(x,T)\in\Rn\times\R_{+}$.
In addition, since $\beta>2$, we learn that
$$D^{2}V(x;T)=\begin{cases}
\frac{\beta(\beta-1)(-x)^{\beta-2}}{\eta}\left(1-e^{-\eta T}\right),&x\le0,\\
0,&x>0,
\end{cases}$$
is continuous. This, combined with the continuity of
$$\partial_{T}V(x;T)=\begin{cases}
(-x)^{\beta}e^{-\eta T},&x\le0,\\
0,&x>0,
\end{cases}$$
ensures that $V\in C^{2,1}(\Rn\times\R_{+},\R)$; thus Assumption \eqref{A1} is satisfied.
For
$$
\phi(x)=\partial_{t}v(t,x;T)|_{t=0}+f(x,u^{c}(x;T))=
\begin{cases}
(-x)^{\beta}\left(1-e^{-\eta T}\right),&x\le0,\\
0,&x>0.
\end{cases}
$$
to be positive when $x<0$ we require $\eta>0$.  Hence Assumption \eqref{A2} is satisfied when $\eta>0$.
In light of the continuity of
$V$ we know that for any $\epsilon>0$ there is $\delta>0$ such that $V(x;T)<\epsilon$ for $-\delta<x<0$;
thus Assumption \eqref{A3} is satisfied. For $x_{1}<x_{2}<0$ we have $V(x_{1};T)>V(x_{2};T)>0$
thus Assumption \eqref{A4}
is satisfied if we choose $h$ to be $V(\cdot;T)$ itself.

In conclusion, for given $r$, $b$ and $\sigma$, we obtain that
Assumptions (\ref{A1}-\ref{A4}) are satisfied for all choices of $\beta$ such that:
\begin{equation}\label{ex:betarange}
2\,<\,\beta\,<\, 1+\frac{1}{2}\frac{1}{r}\frac{(b-r)^2}{\sigma^2},
\end{equation}
where the inequality on the right-hand side corresponds to the condition $\eta>0$.
Thus, a necessary condition to have a stabilizing control process is that the right-hand side
of the above inequality is greater that the left-hand side; that is $(b-r)^2>2r\sigma^2$.
Finally, it is easy to se that the constraint $u\in[c_{1},c_{2}]$ can be met
provided that it is possible to choose
\begin{equation}\label{ex:betaconst}
\beta\,>\, 1+\frac{(b-r)}{|c_1|\sigma^2},
\end{equation}
which, taking into account \eqref{ex:betarange}, can be done if $|c_1|>2r/(b-r)$.

\subsubsection{Terminal cost}\label{sec:exterm}

It is also possible to consider a cost function in the form
\begin{equation}
f(x,u):=0,\,\forall x,
\quad\mbox{ and }\quad
g(x):=\begin{cases}
(-x)^{\beta},&x\le0,\\
0,&x>0,
\end{cases}
\end{equation}
The derived controlled process turns out to be the same as in \eqref{eq:ex1uc}.
However, in this case, the value function is given by
\begin{eqnarray}\label{eq:ex2value}
v(t,x;T)=\begin{cases}
(-x)^{\beta}e^{\eta(t-T)},&x\le0,\\
0,&x>0.
\end{cases}
\end{eqnarray}
where $\eta$ is given again by \eqref{eq:ex1eta}.
Through similar steps, one eventually obtains the same
stability conditions of the previous case.\\

\subsection{Numerical illustration}\label{sec:exsim}

\begin{figure}[t]
\centering
\includegraphics[width=0.78\columnwidth]{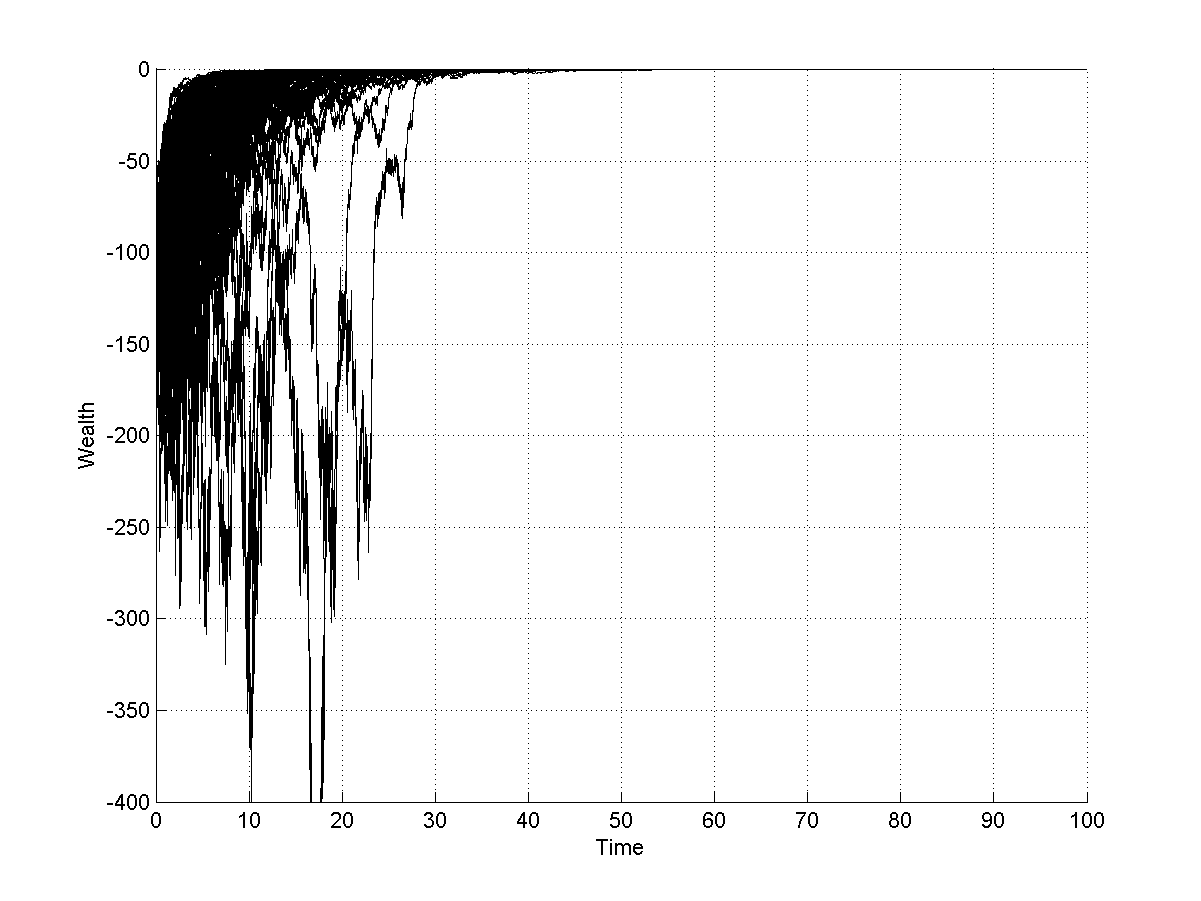}
\caption{Wealth process for $\beta=2.1$ (100 simulations)}
\label{fig:b21}
\end{figure}

\begin{figure}[h!]
\centering
\includegraphics[width=0.78\columnwidth]{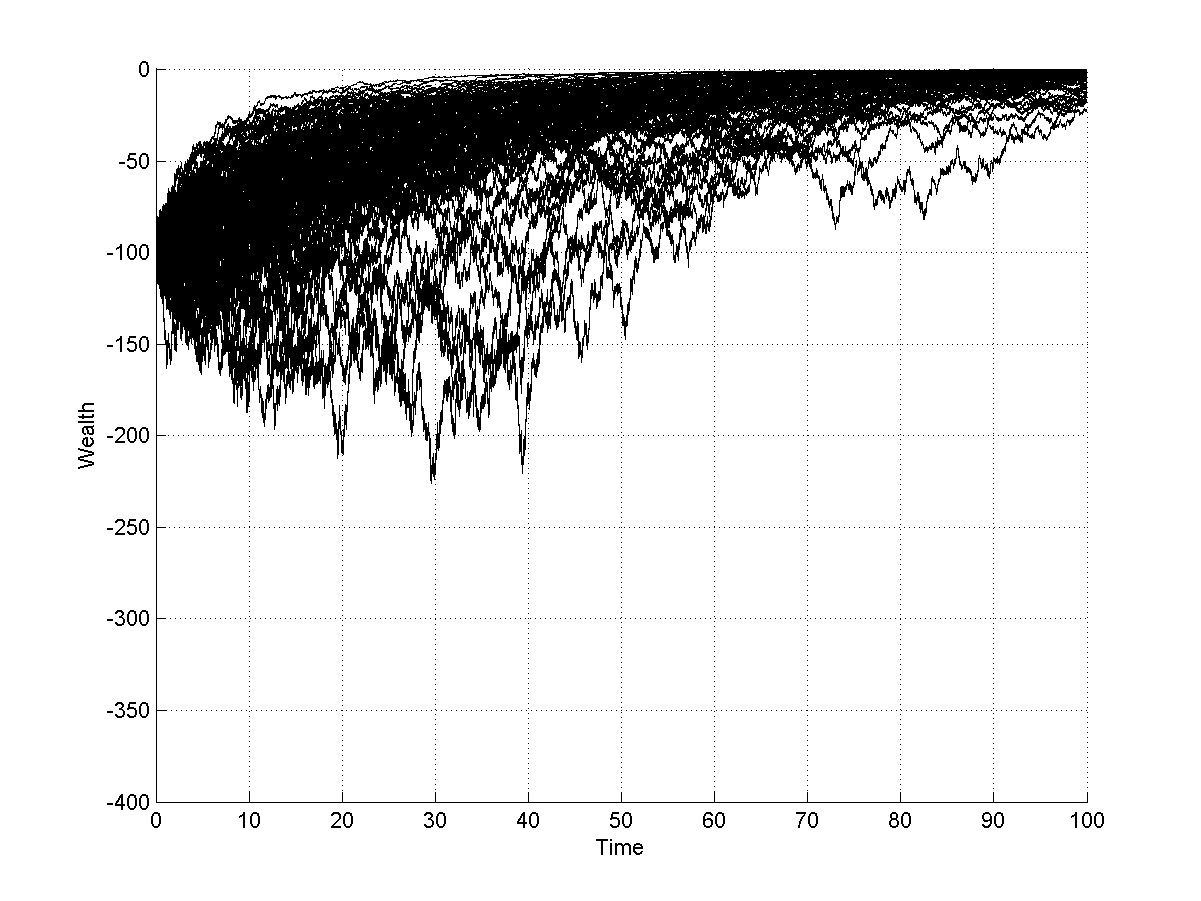}
\caption{Wealth process for $\beta=4.5$ (100 simulations)}
\label{fig:b45}
\end{figure}

\begin{figure}[t]
\centering
\includegraphics[width=0.78\columnwidth]{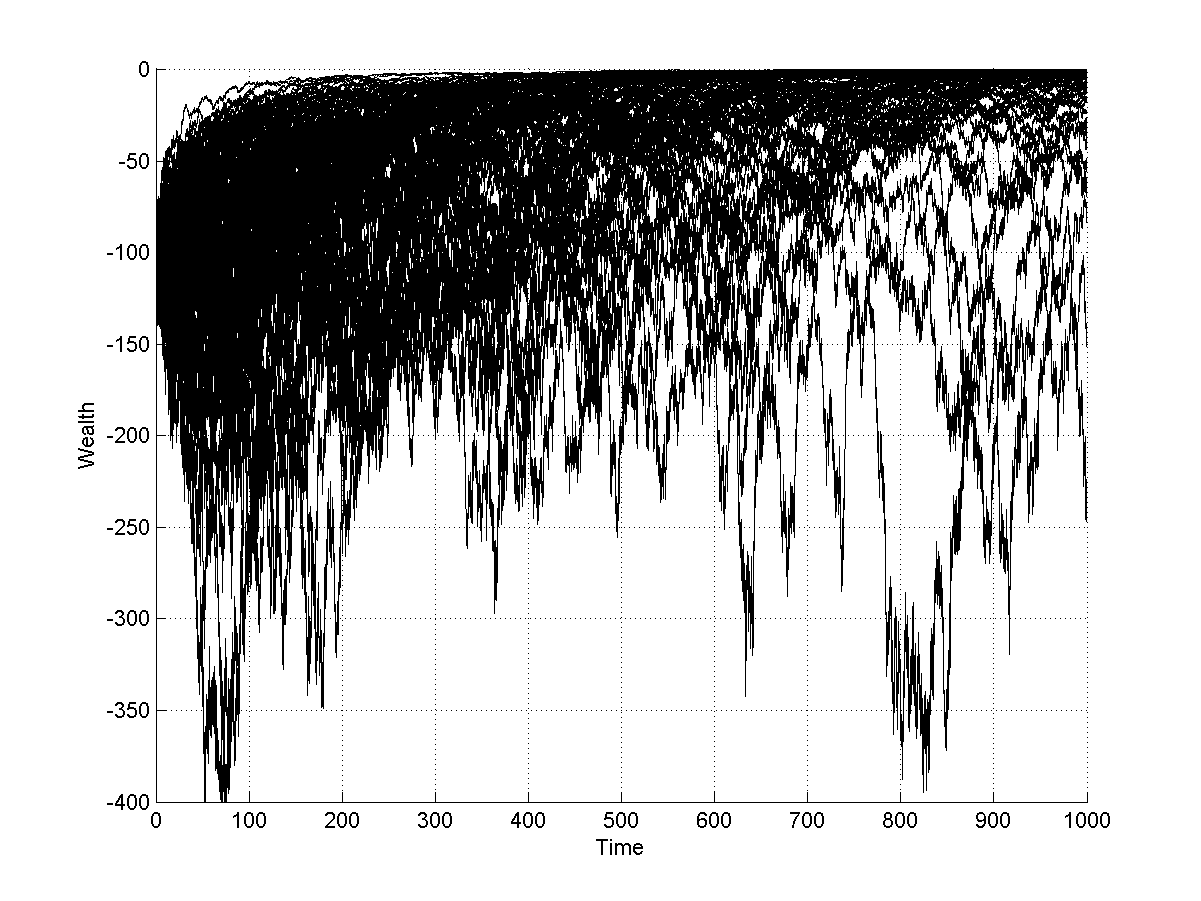}
\caption{Wealth process for $\beta=7.8$ (100 simulations)}
\label{fig:b80}
\end{figure}

We present a simulation example where:
$r=0.03$, $b=0.1$, $\sigma = 0.15$,  and $x_0 = -100$.\linebreak
Here we illustrate the behaviour of the wealth process under the
receding horizon control process \eqref{eq:ex1uc}
for three  different choices of $\beta$ in  cost function \eqref{eq:excost}: $\beta=2.1$, $\beta=4.5$ and $\beta=8$.
The corresponding Monte Carlo simulations are displayed in
 Figures \ref{fig:b21}-\ref{fig:b80} respectively.
Note that, according to \eqref{ex:betarange}, for the given values of $r$, $b$ and $\sigma$,
we have that the  wealth process is asymptotically  stable for $\beta\in(2,\,4.6)$.
By inspecting the figures, it can be seen that for $\beta=2.1$
the wealth  process is clearly asymptotically stable
but there is a significant risk of a large initial undershoot.
For $\beta=4.5$ the process converges much slower but the risk
of a initial undershoot is reduced.
For $\beta=7.8$ the wealth process is not
asymptotically stable (note the different time scale in the figure).

\clearpage


\section{Conclusions}

In this note, we have discussed the RHC strategy for systems
described by continuous-time SDEs.
We have obtained conditions on the associated
finite-horizon optimal control problem
which guarantee
the asymptotic stability of the RHC law.
We have shown that these conditions recall their
deterministic counterpart.
We have illustrated the results with
a simple example in which the RHC law
can be obtained explicitly.
In current work, we are addressing
the problem of implementation in more realistic applications.
For this purpose, it will be necessary to formulate
conditions on the control problem \eqref{eq:cost}
which can guarantee that Assumptions (\ref{A1}-\ref{A4})
hold true and which can be imposed or verified easily.
For problems where the dimension of the
state space is not prohibitive,
it is  possible to solve the associated finite-horizon
optimal control problem with numerical methods \cite{Pham-05,Borkar-05}.
The other possibility is is to approach the
finite-horizon problem \eqref{eq:cost}
by direct on-line optimization.
In this case, one considers a parameterized class
of feedback policies and iteratively optimizes
the parameter of the feedback policy, at regular time intervals,
conditioned on the value of the current state.
In this context simulation-based optimization
methods have shown to be  promising tools, see e.g.~\cite{Maciejowski-et-al-05,Lecchini-et-al-06,Kantas-et-al-10}.\\

\noindent{\bf References}

\bibliographystyle{model1-num-names}
\bibliography{../../../BIB/Biblio}

\end{document}